\documentclass[12pt]{amsart}
\usepackage{amsmath, amssymb, latexsym, amsthm, mathrsfs, color, mathtools, latexsym}
\usepackage[hyphens]{url}

\usepackage[top=2.5cm, bottom=2.5cm, left=2.5cm, right=2.5cm]{geometry}

\newcommand{\nc}{\newcommand}

\nc{\thusfar}{\par\bigskip\centerline{\my{--- Edited thus far ---}}\par\bigskip}
\nc{\lei}{\le^\oo}
\nc{\card}[1]{\left|#1\right|}
\nc{\medcard}[1]{\biggl|\,#1\,\biggr|}
\nc{\smallcard}[1]{|\,#1\,|}
\nc{\bds}{bidirectional $\roth$-scale}
\nc{\bfP}{\mathbf{P}}
\nc{\bfQ}{\mathbf{Q}}
\nc{\bbT}{\mathbb{T}}
\nc{\bbN}{\mathbb{N}}
\nc{\bbC}{\mathbb{C}}
\nc{\beq}{\begin{eqnarray*}}\nc{\eeq}{\end{eqnarray*}}
\nc{\mbq}{\mb{?}}
\nc{\mb}[1]{{\mbox{\textbf{#1}}}}
\nc{\nop}{$\times$}
\nc{\fbn}{\!\!\fbox{\!\nop\!}\!\!}
\nc{\yup}{\checkmark}
\nc{\forces}{\Vdash}
\nc{\name}[1]{\dot{#1}}
\nc{\tf}{\my{FINISHED THUS FAR}}
\nc{\FU}{Fr\'echet--Urysohn}
\nc{\gs}{$\gamma$~space}
\nc{\Ga}{\Gamma}\nc{\Om}{\Omega}
\nc{\smallbinom}[2]{\begin{psmallmatrix} #1\\ #2 \end{psmallmatrix}}
\nc{\bgamma}{\smallbinom{\Om}{\Ga}}
\newcommand{\two}{\{0,1\}}
\nc{\productive}[2]{(#1,\allowbreak #2)^\x}
\nc{\prdct}[1]{(#1)^\x}
\nc{\Sel}{\mathsf{S}}
\nc{\sset}[2]{\{\,#1 : #2\,\}}
\nc{\smb}[1]{{\!\!\mb{#1}\!\!}}
\nc{\medset}[2]{{\biggl\{\,#1 : #2\,\biggr\}}}
\nc{\smallmedset}[2]{{\bigl\{\,#1 : #2\,\bigr\}}}
\nc{\set}[2]{{\left\{\,#1 : #2\,\right\}}}
\nc{\seq}[2]{{\la\, #1 : #2\,\ra}}
\nc{\eseq}[1]{#1_1, \allowbreak #1_2, \allowbreak\dotsc} 
\nc{\cube}{(\Cantor)^\bbN}
\nc{\Match}{\op{Match}}
\nc{\concat}[1]{\hat{\phantom{a}}\langle #1\rangle}
\nc{\poset}{\mathbb{P}}
\nc{\fn}[1]{{\op{Fn}(#1\times\w,2)}}
\nc{\linadd}{\op{linadd}}
\nc{\nonprod}{\non^\x}
\nc{\alephes}{{\aleph_0}}
\nc{\my}[1]{\marginpar{\textcolor{red}{***}}\textcolor{red}{#1}}
\nc{\later}[1]{{\color{green} #1}}
\nc{\BTs}[1]{{\color{green} #1 (BT)}}
\nc{\Cp}{\op{C}_\mathrm{p}}
\nc{\Bp}{\op{B}_p}
\nc{\Pa}[8]{\bibitem{#1} {#2}, \emph{#3}, {#4} \textbf{#5} ({#6}), {#7}--{#8}.}
\nc{\tPa}[5]{\bibitem{#1} {#2}, \emph{#3}, {#4}, to appear.}
\nc{\sPa}[4]{\bibitem{#1} {#2}, \emph{#3}, {#4}, submitted.}
\nc{\Bc}[9]{\bibitem{#1} {#2}, \emph{#3}, in: \textbf{#4} (#5), #6 #7, #8--#9.}
\nc{\fD}{\mathfrak{D}}
\nc{\fX}{\mathfrak{X}}
\nc{\Onbd}{\Op_{\mathrm{nbd}}} 
\nc{\Omnb}{\Om_{\mathrm{nbd}}} 
\nc{\od}{\mathfrak{od}}
\nc{\Setting}[7]{\xymatrix@R=4pt@C=7pt{#1\ar@{-}[r]&#2\ar@{-}[r]&#3\\&#4\ar@{-}[u]\\
#5\ar@{-}[uu]\ar@{-}[r] & #6\ar@{-}[u]\ar@{-}[r] & #7\ar@{-}[uu]}}
\nc{\mx}[1]{\begin{matrix}#1\end{matrix}}
\nc{\plim}{p\txt{-}\lim}
\nc{\Bgp}{{\Z^\bbN}}
\nc{\Cgp}{{{\Z_2}^\bbN}}
\nc{\Cite}[1]{\textbf{[#1]}}
\nc{\Next}[1]{{#1^+}}
\nc{\cFin}{\mathrm{cF}}
\nc{\scsp}{\text{-scale space}}
\nc{\cfn}{\text{cofinal}\ }
\nc{\Con}{\text{Concentrated}}
\nc{\Lind}{\text{Lindel\"of}\,}
\nc{\con}{\text{-Concentrated}}
\nc{\lind}{\text{-Lindel\"of}\,}
\nc{\ctbl}{\text{countably }\allowbreak}
\nc{\Men}{\text{Menger}}
\nc{\men}{\text{-Menger}}
\nc{\Hur}{\text{Hurewicz}}
\nc{\intvl}[2]{{[#1(#2),\allowbreak #1(#2\!+\!1))}}
\nc{\Bdd}{\mathbf{B}}
\nc{\Dfin}{\mathfrak{D}_\mathrm{fin}}
\nc{\grbl}{{\mbox{\textit{\tiny gp}}}}
\nc{\bbP}{\mathbb{P}}
\nc{\BOfat}{\B_{\Om_{\mathrm{fat}}}}
\nc{\Bgood}{\B_{\mathrm{good}}}
\nc{\compactN}{\cl{\mathbb{N}}}
\nc{\blocks}[2]{\op{cl}_{#2}(#1)}
\nc{\blocksplus}[2]{\op{cl}^+_{#2}(#1)}
\nc{\arx}[1]{\texttt{http://arxiv.org/math/#1}}
\nc{\bq}{\begin{quote}}
\nc{\eq}{\end{quote}}
\nc{\cl}[1]{\overline{#1}}
\nc{\CH}{the Continuum Hypothesis}
\nc{\MA}{Martin's Axiom}
\nc{\Bfat}{\B_\mathrm{fat}}
\nc{\inv}{^{-1}}
\nc{\Cantor}{{\two^\bbN}}
\nc{\bP}{\mathbf{P}}
\nc{\bof}{\op{\fb}}
\nc{\dof}{\op{\fd}}
\nc{\bofF}{\bof(\cF)}
\nc{\sr}[3]{\underset{\mbox{#3}}{\mbox{#1}}}
\nc{\gp}{\binom{\Om}{\Ga}}
\nc{\gpsmall}{\mbox{$\gp$}}
\nc{\gig}{\gimel}
\nc{\gns}{\sone(\Om,\gig)}
\nc{\nsr}[2]{#1}
\nc{\Srg}{{\mathbb{S}}}
\nc{\Srgs}{{\mathbb{S}^*}}
\nc{\NN}{{\bbN^{\bbN}}}
\nc{\ZN}{{\Z^{\bbN}}}
\nc{\NNup}{{\bbN^{\uparrow\bbN}}}
\nc{\Pof}{\op{P}}
\nc{\PN}{{\Pof(\bbN)}}
\nc{\roth}{{[\bbN]^{\mbox{\tiny $\infty$}}}} 
\nc{\Fin}{[\bbN]^{\text{$<\!\!\infty$}}} 
\nc{\ici}{[\bbN]^{ \infty, \infty}}
\nc{\Inc}{{\compactN^{\uparrow\bbN}}}
\nc{\powInc}[1]{{\big(\Inc\big)^{#1}}}
\nc{\powFin}[1]{{\big(\Fin\big)^{#1}}}
\nc{\powPN}[1]{{\big(\PN\big)^{#1}}}
\nc{\NcompactN}{{\compactN^\bbN}}
\nc{\Uarrow}{\smash{\big\uparrow}}
\nc{\LE}{\preccurlyeq}
\nc{\GE}{\succcurlyeq}
\nc{\op}{\operatorname}
\nc{\im}{\op{im}}
\nc{\Span}{\op{span}}
\nc{\maxfin}{\op{maxfin}}
\nc{\ran}{\op{range}}
\nc{\iso}{\cong}
\nc{\Madd}{{\M}^\star}
\nc{\cI}{\mathcal{I}}
\nc{\cJ}{\mathcal{J}}
\nc{\scrA}{\mathscr{A}}
\nc{\scrB}{\mathscr{B}}
\nc{\scrC}{\mathscr{C}}
\nc{\scrD}{\mathscr{D}}
\nc{\scrF}{\mathscr{F}}
\nc{\scrK}{\mathscr{K}}
\nc{\A}{\forall}
\nc{\B}{\mathrm{B}}
\nc{\cB}{\mathcal{B}}
\nc{\bB}{\mathbf{B}}
\nc{\BS}{\mathbf{B}(\mathcal{S})}
\nc{\BF}{\mathbf{B}(\mathcal{F})}
\nc{\BU}{\mathbf{B}(\mathcal{U})}
\nc{\cSp}{\mathcal{S}^+}
\nc{\cFp}{\mathcal{F}^+}
\nc{\cUp}{\mathcal{U}^+}

\nc{\BG}{\B_\Ga}
\nc{\BL}{\B_\Lambda}
\nc{\BT}{\B_\Tau}
\nc{\BTstar}{\B_{\Tau^*}}
\nc{\BO}{\B_\Om}
\nc{\DO}{\cD_\Om}
\nc{\KO}{\cK_\Om}
\nc{\CG}{C_\Ga}
\nc{\CL}{C_\Lambda}
\nc{\CT}{C_\Tau}
\nc{\CTstar}{C_{\Tau^*}}
\nc{\CO}{C_\Om}
\nc{\COgp}{C_{\Om^{\grbl}}}
\nc{\CLgp}{C_{\Lambda^{\grbl}}}
\nc{\BOgp}{\B_{\Om}^{\grbl}}
\nc{\BLgp}{\B_{\Lambda^{\grbl}}}
\nc{\sfC}{\mathsf{C}}
\nc{\sfD}{\mathsf{D}}
\nc{\bD}{\mathbf{D}}
\nc{\Tau}{\mathrm{T}}
\nc{\cA}{\mathcal{A}}
\nc{\cK}{\mathcal{K}}
\nc{\cD}{\mathcal{D}}
\nc{\cF}{\mathcal{F}}
\nc{\cS}{\mathcal{S}}
\nc{\cT}{\mathcal{T}}
\nc{\cG}{\mathcal{G}}
\nc{\cY}{\mathcal{Y}}
\nc{\J}{\mathcal{J}}
\nc{\cL}{\mathcal{L}}
\nc{\cM}{\mathcal{M}}
\nc{\cN}{\mathcal{N}}
\nc{\cH}{\mathcal{H}}
\nc{\cO}{\mathcal{O}}
\nc{\Op}{\mathrm{O}}
\nc{\rmA}{\mathrm{A}}
\nc{\rmF}{\mathrm{F}}
\nc{\rmB}{\mathrm{B}}
\nc{\rmD}{\mathrm{D}}
\nc{\rmP}{\mathrm{P}}
\nc{\cC}{\mathcal{C}}
\nc{\cP}{\mathcal{P}}
\nc{\bbQ}{\mathbb{Q}}
\nc{\bbR}{\mathbb{R}}
\nc{\cU}{\mathcal{U}}
\nc{\Un}{\bigcup}
\nc{\cV}{\mathcal{V}}
\nc{\cW}{\mathcal{W}}
\nc{\Z}{{\mathbb Z}}
\nc{\Impl}{\Rightarrow}
\long\def\forget#1\forgotten{\marginpar{\textcolor{green}{Forgetting...}}}
\nc{\ft}{\mathfrak{t}}
\nc{\fb}{\mathfrak{b}}
\nc{\fc}{\mathfrak{c}}
\nc{\fd}{\mathfrak{d}}
\nc{\fg}{\mathfrak{g}}
\nc{\oo}{\infty}
\nc{\fr}{\mathfrak{r}}
\nc{\fk}{\mathfrak{k}}
\nc{\bidi}{\mathfrak{bidi}}
\nc{\fu}{\mathfrak{u}}
\nc{\fh}{\mathfrak{h}}
\nc{\fp}{\mathfrak{p}}
\nc{\fj}{\mathfrak{j}}
\nc{\fs}{\mathfrak{s}}
\nc{\w}{\omega}
\nc{\x}{\times}
\nc{\Iff}{\Leftrightarrow}
\newcommand\comp{^{\text{\tt c}}}
\nc{\nin}{\notin}
\nc{\cat}{\hat{\ }}
\nc{\sub}{\subseteq}
\nc{\spst}{\supseteq}
\nc{\sm}{\setminus}
\nc{\as}{\subseteq^*}
\nc{\les}{\le^*}
\nc{\leinf}{\le^{\infty}}
\nc{\leS}{\le_S}
\nc{\leF}{\le_F}
\nc{\leU}{\le_U}
\nc{\rest}{\restriction}

\nc{\la}{\langle}
\nc{\ra}{\rangle}
\nc{\E}{\exists}
\nc{\dom}{\op{dom}}
\nc{\cov}{\op{cov}}
\nc{\add}{\op{add}}
\nc{\cof}{\op{cof}}
\nc{\cf}{\op{cf}}
\nc{\non}{\op{non}}
\nc{\unif}{\op{non}}
\nc{\COV}{\op{COV}}
\nc{\ADD}{\op{ADD}}
\nc{\COF}{\op{COF}}
\nc{\NON}{\op{NON}}
\nc{\impl}{\to}
\nc{\Lp}{\mathcal{L_\p}}
\nc{\Wlog}{without loss of generality}
\newtheorem{thm}{Theorem}[section]
\nc{\bthm}{\begin{thm}} \nc{\ethm}{\end{thm}}
\newtheorem{prop}[thm]{Proposition}
\nc{\bprp}{\begin{prop}} \nc{\eprp}{\end{prop}}
\newtheorem{fact}[thm]{Fact}
\nc{\bfct}{\begin{fact}} \nc{\efct}{\end{fact}}
\newtheorem{prob}[thm]{Problem}
\nc{\bprb}{\begin{prob}} \nc{\eprb}{\end{prob}}
\newtheorem{lem}[thm]{Lemma}
\nc{\blem}{\begin{lem}} \nc{\elem}{\end{lem}}
\newtheorem{app}[thm]{Application}
\nc{\bapp}{\begin{app}} \nc{\eapp}{\end{app}}
\newtheorem{claim}[thm]{Claim}
\nc{\bclm}{\begin{claim}} \nc{\eclm}{\end{claim}}
\newtheorem{cor}[thm]{Corollary}
\nc{\bcor}{\begin{cor}} \nc{\ecor}{\end{cor}}
\newtheorem{conj}[thm]{Conjecture}
\nc{\bcnj}{\begin{conj}} \nc{\ecnj}{\end{conj}}
\theoremstyle{definition}
\newtheorem{defn}[thm]{Definition}
\nc{\bdfn}{\begin{defn}} \nc{\edfn}{\end{defn}}
\newtheorem{obs}[thm]{Observation}
\nc{\bobs}{\begin{obs}} \nc{\eobs}{\end{obs}}
\theoremstyle{remark}
\newtheorem{rem}[thm]{Remark}
\nc{\brem}{\begin{rem}} \nc{\erem}{\end{rem}}
\newtheorem{cnv}[thm]{Convention}
\nc{\bcnv}{\begin{cnv}} \nc{\ecnv}{\end{cnv}}
\newtheorem{exam}[thm]{Example}
\nc{\bexm}{\begin{exam}} \nc{\eexm}{\end{exam}}
\nc{\bpf}{\begin{proof}} \nc{\epf}{\end{proof}}
\nc{\be}{\begin{enumerate}}
\nc{\ee}{\end{enumerate}}
\nc{\bi}{\begin{itemize}}
\nc{\bimy}{\my{\begin{itemize}}
\nc{\eimy}{\end{itemize}}}
\nc{\itm}{\item}
\nc{\ei}{\end{itemize}}
\nc{\Subsection}[1]{\goodbreak\subsection*{#1}}

\nc{\sone}{\mathsf{S}_1}
\nc{\sfin}{\mathsf{S}_\mathrm{fin}}
\nc{\ufin}{\mathsf{U}_\mathrm{fin}}
\nc{\Split}{\mathsf{Split}}

\nc{\gone}{\mathsf{G}_1}    \nc{\gfin}{\mathsf{G}_\mathrm{fin}}

\title[Products of general Menger spaces]{Products of general Menger spaces}

\subjclass[2010]{
54D20, 03E17.
}

\author[P. Szewczak]{Piotr Szewczak}
\address{Piotr Szewczak, 
	Institute of Mathematics, Faculty of Mathematics and Natural Science College of Sciences, Cardinal Stefan Wyszy\'nski University in Warsaw, Poland,
	and
	Department of Mathematics, Bar-Ilan University, Israel
}
\email{p.szewczak@wp.pl}
\urladdr{http://piotrszewczak.pl}

\author[B. Tsaban]{Boaz Tsaban}
\address{Boaz Tsaban,
	Department of Mathematics, Bar-Ilan University, Israel}
\email{tsaban@math.biu.ac.il}
\urladdr{http://math.biu.ac.il/~tsaban}

\begin{document}

\begin{abstract}
We study products of general topological spaces with Menger's covering property, 
and its refinements based on filters and semifilters.
To this end, we extend the projection method from
the classic real line topology to the Michael topology.
Among other results, we prove that, assuming \CH{}, 
every productively Lindel\"of space is productively Menger,
and every productively Menger space is productively Hurewicz.
None of these implications is reversible.
\end{abstract}

\maketitle

\section{Introduction}

Let $\roth$ be the set of all infinite subsets of the set $\bbN$ of natural numbers. For sets $a, b\in\roth$, we write $a\as b$ if the set $a\sm b$ is finite.
A \emph{semifilter}~\cite{BZCSFQdM} is a set $S\sub\roth$ such that, for each set
$s\in S$ and each set $b\in\roth$ with $s\as b$, we have $b\in S$.  For a semifilter $S$, let $S^+:=\sset{a\in\roth}{a\comp\notin S}$.  
For all sets $a\in S$ and $b\in S^+$, the intersection  $a\cap b$ is infinite.

By \emph{space} we mean a topological space.
Let $S$ be a semifilter. 
A  space $X$
 is \emph{$S$-Menger}~\cite{sfh, pMReal} if for each sequence $\eseq{\cU}$ of open covers of the space $X$, there are finite sets $\cF_1\sub\cU_1$, $\cF_2\sub\cU_2$, \dots such that $\sset{n\in\bbN}{x\in\Un\cF_n}\in S$ for all points $x\in X$. 
A  space is \emph{Menger}~\cite{Menger24, Hure25} if it is $\roth$-Menger~\cite{Hure25}. 
Let $\cFin:=(\roth)^+$, the (semi)filter of cofinite subsets of $\bbN$. A space is \emph{Hurewicz}~\cite{Hure25} if it is $\cFin$-Menger. 
Restricting the definition of $S$-Menger spaces to countable open covers, we obtain the definition of \emph{countably $S$-Menger} space.
This makes it possible to extend our investigations beyond the Lindel\"of realm.

\bprp\label{prp:LindctblSMen=SMen}
Let $S$ be a semifilter. 
A space $X$ is $S$-Menger if and only if it is Lindel\"of and countably $S$-Menger.\qed
\eprp

We identify each set $a\in\roth$ with its increasing enumeration.
Thus, for a natural number~$n$, $a(n)$ is the $n$-th smallest element of the set $a$. 
We identify the Cantor space $\Cantor$ with the family $\PN$ of all subsets of the set $\bbN$.
We denote by $\Fin$ the family of finite subsets of $\bbN$, so that $\PN=\roth\cup\Fin$.
Since we identify every set $a\in\roth$ with its increasing enumeration, we have  $\roth\sub\NN$. 
The topology of the space $\roth$ (a subspace of the Cantor space $\PN$) coincides
with the subspace topology induced by the Baire space~$\NN$. 

Let $S$ be a semifilter. 
For elements $a,b\in\roth$, we write $a\leS b$ if $\sset{n\in\bbN}{a(n)\le b(n)}\in S$.
Then $a\nleq_S b$ if and only if $b<_{S^+} a$. 
We write $a\les b$ if $a\le_\cFin b$, $a\leinf b$ if $a\le_\roth b$, and $a\le b$ if 
$a(n)\leq b(n)$ for all natural numbers $n$.

A map $\Psi$ from a space $X$ into $\roth$ is \emph{upper continuous} if the sets
$\sset{x\in X}{\Psi(x)(n)\le m}$ are open for all natural numbers $n$ and $m$.
In particular, continuous functions are upper continuous. 

Let $X$ and $Y$ be spaces. A set-valued map $\Phi\colon X\Rightarrow Y$ is 
\emph{compact-valued upper semicontinuous (cusco)} if, 
for each point $x\in X$, the set $\Phi(x)$ is a nonempty compact subset of the space $Y$, 
and for every open set $V\sub Y$, the set $\Phi\inv[V]:=\sset{x\in X}{\Phi(x)\sub V}$ 
is open in the space~$X$. 
The \emph{image} of a set~$X$ under a set-valued map $\Phi$ is the set $\Phi[X]:=\Un_{x\in X}\Phi(x)$.
The following observation can be proved using earlier methods~\cite[Theorem~7.3]{SPMProd}.

\bprp\label{prp:ctblSMenChar}
Let $X$ be a  space, and $S$ be a semifilter.
The following assertions are equivalent:
\be
\itm The space $X$ is countably $S$-Menger.
\itm Every upper continuous image of the space $X$ in $\roth$
is $\leS$-bounded.
\itm Every cusco image of the space $X$ in $\roth$
is $\leS$-bounded.\qed
\ee
\eprp

For a semifilter $S$, 
let $\bof(S)$ be the minimal cardinality of a $\leS$-unbounded subset of $\roth$. 

\bcor\mbox{}
\be
\itm Every space of cardinality smaller than $\bof(S)$ is countably $S$-Menger. 
\itm The discrete space of cardinality $\bof(S)$ is not countably $S$-Menger. \qed
\ee
\ecor

\section{Countable covers}

Let $S$ be a semifilter. 
A set $X\sub\roth$ with $\card{X}\ge\bof(S)$
is an \emph{$S$-scale}~\cite[Definition~4.1]{pMReal} if, for each element $b\in\roth$, there is
an element $c\in\roth$ such that $b\le_{S^+} c$ and $c\le_S x$ for all but less than $\bof(S)$ elements $x\in X$. For every semifilter $S$,  $S$-scales provably exist~\cite[Lemma~2.9]{sfh}.
A set $X\sub\roth$ with $\card{X}\ge\bof(S)$ is a \emph{cofinal $S$-scale}~\cite[Definition~6.1]{pMReal} if for each element $b\in\roth$, $b\le_S x$ for all but less than $\bof(S)$ elements $x\in X$. 
By \emph{filter} we mean a semifilter that is closed under finite intersections. If $F$ is a filter, then there is a cofinal $F$-scale if and only if $\bof(F)=\bof(F^+)$~\cite[Corollary~6.3]{pMReal}.

Let $\bfP$ and $\bfQ$ be properties of spaces.
A space $X$ satisfies $\productive{\bfP}{\bfQ}$ if for each space $Y$ with the property $\bfP$, the product space $X\x Y$ has the property $\bfQ$.
Define $\bfP^\x:=\productive{\bfP}{\bfP}$. A space is \emph{productively $\bfP$} if it satisfies $\bfP^\x$.

\subsection{General semifilters}

Let $\kappa$ be an uncountable cardinal number. 
A space $X$ is \emph{$\kappa$-Lindel\"of} if every open cover of $X$ has a subcover of cardinality smaller than $\kappa$. 
A space $X$ with $\card{X}\ge\kappa$ is \emph{$\kappa$-concentrated} on a set $D\sub X$ if $\card{X\sm U}<\kappa$ for all open sets $U$ containing~$D$.

\bthm\label{thm:ctblS}
Let $S$ be a semifilter with a cofinal $S$-scale.
\[
\productive{\bof(S)\con}{\bof(S)\lind}\sub \ctbl S\men.
\]	
\ethm 

For the proof of Theorem~\ref{thm:ctblS}, and for  later discussions, we introduce several notions and auxiliary results. 
Let $X\sub\roth$.
The \emph{Michael topology}~\cite{michael} on the set $X\cup\Fin\sub\PN$ is the one where the points of the set $X$ are isolated, and the neighborhoods of the points of the set $\Fin$ are those induced by the Cantor space topology on $\PN$. 

Let $\kappa$ be an uncountable cardinal number. A set $X\sub\roth$ with $\card{X}\ge\kappa$ is \emph{$\kappa$-unbounded} if $\card{X}\ge\kappa$, and the cardinality of every $\leq$-bounded subset of the set $X$ is smaller than $\kappa$.
A standard argument~\cite[Lemma~2.3]{pMReal} implies the following result.

\blem\label{lem:MichaelConc}
Let $\kappa$ be an uncountable cardinal number, and $X$ be a subset of $\roth$.
The set $X$ is $\kappa$-unbounded if and only if the space $X\cup\Fin$, with the Michael topology,
is $\kappa$-concentrated on the set $\Fin$.\qed
\elem

\bdfn
A \emph{real space} is a subspace of the Cantor space $\PN$.
\edfn

\bprp\label{prp:not(k-con,k-lind)}
Let $\kappa$ be an uncountable cardinal number, and $S$ be a semifilter. Let $X\sub\roth$ be a real space containing a $\kappa$-unbounded set $Y$. 
For the Michael topology on the set $Y\cup\Fin$, the product space $X\x(Y\cup\Fin)$ is not $\kappa$-Lindel\"of.
\eprp
\bpf
Let $X$ be a real space containing a $\kappa$-unbounded set $Y$. By Lemma~\ref{lem:MichaelConc}, 
the space $Y\cup\Fin$ with the Michael 
topology is $\kappa$-concentrated. The diagonal set $\sset{(y,y)}{y\in Y}$ is a closed discrete subset of the product 
space $X\x (Y\cup\Fin)$, of cardinality at least $\kappa$. Thus, the space $X\x (Y\cup\Fin)$ is not 
$\kappa$-Lindel\"of.
\epf

\blem\label{lem:bofSCSunbdd}
Let $S$ be a semifilter with a cofinal $S$-scale. Every $\le_S$-unbounded set in $\roth$ contains a $\bof(S)$-unbounded set.
\elem
\bpf
Let $Y$ be a $\le_S$-unbounded set in $\roth$, and $\sset{d_\alpha}{\alpha<\bof(S)}$ be a cofinal $S$-scale in~$\roth$.
For each ordinal number $\alpha<\bof(S)$, let $y$ be a $\le_S$-bound of
the set $\sset{y_\alpha}{\alpha<\bof(S)}$.
There is an element $y_\alpha\in Y$ such that $y,d_\alpha<_{S^+} y_\alpha$. 
It follows that $y_\alpha\nin\sset{y_\beta}{\beta<\alpha}$.
The set $\sset{y_\alpha}{\alpha<\bof(S)}$ is $\bof(S)$-unbounded: 
Let $b\in\roth$. 
For each ordinal number $\alpha<\bof(S)$ such that $b\le_S d_\alpha$, since $d_\alpha<_{S^+} y_\alpha$, we have $b<^\infty y_\alpha$.  
As the set $\sset{d_\alpha}{\alpha<\bof(S)}$ is 
a cofinal $S$-scale, we have $b<^\oo y_\alpha$
for all but less than $\bof(S)$ ordinal numbers $\alpha<\bof(S)$.
\epf

\begin{proof}[Proof of Theorem~\ref{thm:ctblS}]
Assume that a space $X$ is not countably $S$-Menger. 
By Proposition~\ref{prp:ctblSMenChar}, there is  in $\roth$
a  $\le_S$-unbounded cusco image $Y$ of $X$.
By Lemma~\ref{lem:bofSCSunbdd}, the set $Y$ contains a $\bof(S)$-unbounded set   $Z$. By Proposition~\ref{prp:not(k-con,k-lind)}, 
for the Michael topology on the set $Z\cup \Fin$,
the product space  $Y\x(Z\cup \Fin)$ is not $\bof(S)$-Lindel\"of.  
Since $Y$ is a cusco image of the space $X$, the product space  $X\x(Z\cup\Fin)$ is not $\bof(S)$-Lindel\"of.
\epf

\bdfn
Let $S$ be a semifilter. 
A space $X\cup\Fin$ is a \emph{(cofinal) $S$-scale space} if the set
$X$ is a (cofinal) $S$-scale in $\roth$, and $X\cup\Fin$ has the Michael topology.
\edfn

\brem\label{rem:coarser}
Since moving to a coarser topology on a space provides a continuous image of that space,
all results for (cofinal) $S$-scale spaces are also true when using coarser topologies,
for example, the Cantor space topology.
\erem

A semifilter $S\sub\roth$ is \emph{meager} if $S$ is a meager subset of $\roth$. 
The proof of the next theorem is a  literal repetition of proofs of earlier, analogous results~\cite[Theorems~5.4 and~6.5]{pMReal}.

\bprp\label{prp:ctblProdS}
\mbox{}
\be
\itm Let $S$ be a meager semifilter. Every $S$-scale space is
productively countably Hurewicz.
\itm Every cofinal $\cFin$-scale space is productively countably $S$-Menger,
for all semifilters~$S$.\qed
\ee
\eprp

Every product of a $\fd$-concentrated real space and a Hurewicz space
that is hereditarily Lindel\"of is Menger~\cite[Theorem~4.6]{AddGN}.
The following theorem generalizes that.
For a semifilter $S$, let $\chi(S)$ be the minimal cardinality of a basis for $S$.

\bthm\label{thm:d-concxH=M}
Let $S$ be a semifilter with $\chi(S)<\cf(\bof(S^+))$.
Let $X$ be a Tychonoff $\bof(S^+)$-concentrated space, and $Y$ be an $S$-Menger space. 
Then the product space $X\x Y$ is countably Menger.
\ethm
\bpf
We use the projection method, introduced in an earlier work~\cite{AddGN}.
Let $C$ be a compactification of the space $X$, and 
$Q$ be a countable subset of $X$ on which the space $X$ is concentrated.
Let $\cU_1,\cU_2,\dots$ be countable open covers of the product space $X\x Y$,
by sets open in the product space $C\x Y$.
The product space $Q\x Y$ is a countable union of Menger spaces, and thus Menger.
Let $\cF_1\sub\cU_1,\cF_2\sub\cU_2,\dots$ be finite sets such that the set $U:=\Un_n\Un\cF_n$ contains the set $Q\x Y$. 
Since the set $(C\x Y)\sm U$ is closed in the space $C\x Y$, it is $S$-Menger.
Thus, the projection $H$ of the set $(C\x Y)\sm U$ in $C$ is an $S$-Menger set disjoint of the set $Q$. 

\blem\label{lem:Gdelta}
Let $S$ be a semifilter.
Let $X$ be an $S$-Menger subset of some space, and $G$ a G$_\delta$ set containing $X$.
There are closed sets $C_\alpha$, for the ordinal numbers $\alpha<\chi(S)$, 
such that $X\sub \Un_{\alpha<\chi(S)}C_\alpha\sub G$.
\elem

\bpf 
Let $G:=\bigcap_n G_n$, where every set $G_n$ is open in $C$, and $G_{n+1}\sub G_n$ for all natural numbers $n$. 
Fix a base of cardinality $\chi(S)$ for the semifilter $S$. 
Let $n$ be a natural number. 
For each point $x\in X$, there is an open neighborhood $U^n_x$ of $x$ such that $\cl{U^n_x}\sub G_n$. The family $\cU_n:=\sset{U^n_x}{x\in X}$ is an open cover
of $X$.

Since the space $X$ is $S$-Menger, there are finite sets $\cF_1\sub\cU_1$, $\cF_2\sub\cU_2$, \dots
such that 
\[
I_x:=\smallmedset{n\in\bbN}{x\in \Un\cF_n}\in S
\]
for all points $x\in X$. 

Fix a point $x\in X$. Let $J$ be a member of the basis of $S$ such that $J\sub I_x$.
Let 
\[
K:=\bigcap_{n\in J}\cl{\Un\cF_n}.
\]
Since the set $J$ is infinite and $G_{n+1}\sub G_n$ for all natural numbers $n$, we have $K\sub G$. 
The union of these sets $K$ is as required.
\epf
By Lemma~\ref{lem:Gdelta}, there are closed sets $C_\alpha$ such that $H\sub \Un_{\alpha<\chi(F)}C_\alpha\sub C\sm Q$.
Since the set $X$ is $\bof(S^+)$-concentrated on $Q$, we have $\card{X\cap C_\alpha}<\bof(S^+)$ for all ordinal numbers $\alpha<\chi(S)$, and since $\chi(S)<\cf(\bof(S^+))$, we have 
\[
\medcard{X\cap\Un_{\alpha<\chi(F)} C_\alpha}=
\medcard{\Un_{\alpha<\chi(F)}X\cap C_\alpha}<\bof(S^+).
\]
By Proposition~\ref{prp:ctblSMenChar}, the union of less than $\bof(S^+)$ spaces that are 
countably $S$-Menger is countably Menger.
Thus, the product space $(X\cap \Un_{\alpha<\chi(F)} C_\alpha)\x Y$ is countably Menger, and 
there are finite sets $\cF_1'\sub \cU_1,\cF'_2\sub\cU_2,\dots$ such that the family 
$\Un_n\cF_n'$ covers this product space.
Then the family $\Un_n\cF_n\cup\cF'_n$ is an open cover of the product space $X\x Y$. 
\epf

\subsection{Superfilters}

A \emph{superfilter} is a semifilter $S$ such that, for all sets $a,b\in\roth$, $a\cup b\in S$ implies $a\in S$ or $b\in S$. A semifilter $S$ is a superfilter if
and only if $S=F^+$ for some filter~$F$. 

For a superfilter $S$, every $S$-scale is a cofinal $S$-scale. In particular, every superfilter $S$ has a cofinal $S$-scale~\cite[Proposition~6.6]{pMReal}. 

\bthm\label{thm:ctblSup}
Let $S$ be a superfilter. 
\[
\productive{\text{(cofinal) } S\scsp}{\bof(S)\lind}\sub \ctbl S\men.
\]
\ethm
\bpf 
For superfilters, we have the following improvement of Lemma~\ref{lem:bofSCSunbdd}.

\blem\label{lem:ctblSup}
Let $S$ be a superfilter.
Every $\leS$-unbounded set in $\roth$ contains a (necessarily, cofinal) $S$-scale.
\elem
\bpf
The proof is similar to that of Lemma~\ref{lem:bofSCSunbdd}.
For a superfilter $S$, the intersection of each element of $S$ and each element of $S^+$ belongs to the set $S$.
Thus, if $b\le_S d_\alpha<_{S^+} y_\alpha$, then $b<_{S}y_\alpha$. 
\epf

The rest of the proof is similar to that of Theorem~\ref{thm:ctblS}, 
using Lemma~\ref{lem:ctblSup}.
\epf

Both Theorems~\ref{thm:ctblS} and~\ref{thm:ctblSup} 
imply the following result.

\bcor\label{cor:ctblF+}
Let $S$ be a superfilter.
\[
\productive{\bof(S)\con}{\bof(S)\lind}\sub\ctbl S\men.\qed
\]
\ecor

\subsection{Filters}
Recall that a filter $F$ has a cofinal $F$-scale if and only if
$\bof(F)=\bof(F^+)$.

\bthm\label{thm:ctblF}
Let $F$ be a filter. 
\be
\itm $\productive{F\scsp}{\bof(F)\lind}\sub \ctbl F^+\men$.
\itm If $\bof(F)=\bof(F^+)$, then:
\be
\itm $\productive{\text{cofinal } F\scsp}{\bof(F)\lind} \sub \ctbl F^+\men$.
\itm $\productive{\text{(cofinal) } F^+\scsp}{\bof(F^+)\lind} \sub \ctbl F\men$.
\ee
\ee
\ethm 
\bpf
The following lemma constitutes the combinatorial skeleton of the theorem.

\blem\label{lem:ctblF}
Let $F$ be a filter.
\be
\itm Every $\le_{F^+}$-unbounded set in $\roth$ contains an $F$-scale.
\itm 
If $\bof(F)=\bof(F^+)$, then:
\be
\itm Every $\le_{F^+}$-unbounded set contains a cofinal $F$-scale.
\itm Every $\le_{F}$-unbounded set contains a (cofinal) $F^+$-scale.
\ee
\ee
\elem
\bpf
(1) Let $Y$ be a $\le_{F^+}$-unbounded set in $\roth$, and $\sset{b_\alpha}{\alpha<\bof(F)}$ be a $\le_{F}$-unbounded set in $\roth$.
For each ordinal number $\alpha<\bof(F)$, there is an element $c\in\roth$ such that
$\sset{b_\beta,y_\beta}{\beta<\alpha}\le_{F} c$. Pick an element $y_\alpha\in Y$ such that
$c<_{F} y_\alpha$. Since $F$ is a filter, we have
\[
\set{b_\beta,y_\beta}{\beta<\alpha}<_{F} y_\alpha.
\]
The set $\sset{y_\alpha}{\alpha<\bof(F)}$ is an $F$-scale:
Let $b\in\roth$. As the set $\sset{b_\alpha}{\alpha<\bof(F)}$ is 
$\le_{F}$-unbounded, there is an ordinal number $\alpha<\bof(F)$ such that 
$b<_{F^+} b_\alpha$.
For each ordinal number $\beta<\bof(F)$ that is greater than $\alpha$, we have
$b<_{F^+} b_\alpha<_F y_\beta$.

(2)(a) 
Let $Y$ be a $\le_{F^+}$-unbounded set in $\roth$, and $\sset{d_\alpha}{\alpha<\bof(F)}$ be a cofinal $F$-scale in~$\roth$.
For each ordinal number $\alpha<\bof(F)$, let $y_\alpha$ be an element of $Y$
with $d_\alpha<_F y_\alpha$.
The set $\sset{y_\alpha}{\alpha<\bof(F)}$ is a cofinal $F$-scale: 
Let $b\in\roth$. 
For each ordinal number $\alpha<\bof(F)$ such that $b\le_F d_\alpha$, 
we have  $b\le_F d_\alpha<_F y_\alpha$, and since $F$ is a filter,
$b <_F y_\alpha$. 
The set $\sset{y_\alpha}{\alpha<\bof(F)}$ is $\le_F$-unbounded, and thus its cardinality is not smaller than
$\bof(F)$.

(2)(b)
The proof is similar to that of Lemma~\ref{lem:bofSCSunbdd}. 
Here, since $F$ is a filter, 
$b\leF d_\alpha<_{F^+} y_\alpha$ implies $b<_{F^+}y_\alpha$.
\epf
The proof of the theorem is now similar to that of Theorem~\ref{thm:ctblS},
using Lemma~\ref{lem:ctblF}.
\epf
 
Following the proof of earlier, analogous results~\cite[Theorem~5.3, Theorem~6.2(2)]{pMReal}, 
we obtain the following result.

\bprp\label{prp:ctblProdF}
Let $F$ be a filter.  
\be
\itm  Every $F$-scale space satisfies $\productive{\ctbl F\men}{\ctbl F^+\men}$.
\itm Every cofinal $F$-scale space is productively countably $F$-Menger.\qed
\ee
\eprp

\bthm\label{thm:ctblF2}
Let $F$ be a filter. 
\be 
\itm $\productive{\bof(F)\con}{\bof(F)\lind}\sub\ctbl F^+\men$.	
\itm 	$\productive{\ctbl F^+\men}{\bof(F)\lind}\sub\ctbl F^+\men$.
\itm If $\bof(F)=\bof(F^+)$, then 
\[
\productive{\ctbl F^+\men}{\bof(F^+)\lind}\sub\ctbl F\men.
\]
\ee
\ethm
\bpf
(1) Every $F$-scale is $\bof(F)$-unbounded~\cite[Proposition~4.3]{pMReal}. Apply Lemma~\ref{lem:MichaelConc} and Theorem~\ref{thm:ctblF}(1). 

(2) By Proposition~\ref{prp:ctblProdF}(1), every $F$-scale space is countably $F^+$-Menger. Apply Theorem~\ref{thm:ctblF}(1).

(3) Apply Theorem~\ref{thm:ctblF}(2).
\epf

We only consider \emph{nonprincipal} ultrafilters.
An ultrafilter is simultaneously a filter and a superfilter.

\bcor\label{cor:ctblU} Let $U$ be an ultrafilter. 
\be
\itm $\productive{\bof(U)\con}{\bof(U)\lind}\sub\ctbl U\men$.
\itm $\productive{\ctbl U\men}{\bof(U)\lind}\sub\ctbl U\men$.
\itm  Every $U$-scale space is productively countably $U$-Menger.	
\ee
\ecor
\bpf
(1) Apply any of Theorem~\ref{thm:ctblS}, Theorem~\ref{thm:ctblF}, Theorem~\ref{thm:ctblF2}(1), or Corollary~\ref{cor:ctblF+}.

(2) Apply Theorem~\ref{thm:ctblF2}(2) or~(3).

(3) Apply $U^+=U$ and Proposition~\ref{prp:ctblProdF}(1) or~(2).
\epf

Let $F$ be a filter. By Lemmata~\ref{lem:ctblSup} and~~\ref{lem:ctblF},
every $\le_{F^+}$-unbounded set contains a (cofinal) $F^+$-scale
and an $F$-scale.
However, in general, a $\leF$-unbounded set may not contain an $F$-scale;
this may be the case for the filter $\cFin$, even if there is a cofinal 
$\cFin$-scale, as we now explain.

A \emph{Luzin set} is an uncountable real space 
whose intersection with each meager set is countable. 
In particular, every Luzin set is $\le_{\cFin}$-unbounded. 

\bprp\label{prp:Lu}
No Luzin subset of $\roth$ contains a $\cFin$-scale.
\eprp

Proposition~\ref{prp:Lu} follows from the following Lemma.

\blem
Every $\cFin$-scale is meager in $\roth$.
\elem
\bpf
Let $Y\sub\roth$ be a $\cFin$-scale. Fix a coinfinite set $b\in\roth$. 
There is a set $c\in\roth$ such that $b\leinf c$ and $c\les y$ for all but less than $\fb$ elements $y\in Y$. 
The set $Z:=\sset{y\in Y}{c\les y}$ is a countable union of sets
of the form $\sset{y\in Y}{c'\leq y}$, with $c'\in\roth$ a coinfinite set.

Fix a coinfinite set $c'\in\roth$. The set $X:=\sset{x\in\roth}{c'\leq x}$ is nowhere dense: Let $U$ be an open subset of $\roth$.
Let $a\in U$ be a cofinite set. There is a natural number $n$ such that $a(n)<c'(n)$. The open set 
\[
\sset{x\in\roth}{x(1)=a(1),\dotsc,x(n)=a(n)}
\]
is a subset of the set $U$ disjoint of the set $X$.

Since $\card{Y\sm Z}<\fb$, the remainder $Y\sm Z$ is $\les$-bounded, 
and thus meager.
\epf

\brem
Let $V$ be a model of set theory satisfying \CH{}. In $V$,
let $\kappa$ be a cardinal number with uncountable cofinality,
and $\bbC_\kappa$ be the Cohen forcing notion adding $\kappa$ Cohen reals.
In the extended model $V^{\bbC_\kappa}$, we have 
\[
\aleph_1=\fb=\bof(\cFin) \le \fd=\bof(\cFin^+)=\fc=\kappa,
\]
and the canonical set of generic Cohen reals is a Luzin set of cardinality $\kappa$~\cite[Lemma~8.2.6]{BarJu}.
By Proposition~\ref{prp:Lu},
this shows that, for the filter $F:=\cFin$,
neither $\bof(F)=\bof(F^+)$ (taking $\kappa=\aleph_1$)
nor 
$\bof(F)<\bof(F^+)$ (taking $\kappa=\aleph_2$)
imply that 
every $\leF$-unbounded set in $\roth$ contains an $F$-scale.
\erem

Let $F$ be a filter.
By Lemma~\ref{lem:bofSCSunbdd}, if $\bof(F)=\bof(F^+)$, then 
every $\le_F$-unbounded set contains a $\bof(F)$-unbounded set.
The assumption $\bof(F)=\bof(F^+)$ is not redundant. Here too, 
a consistent counterexample exists for the filter $F:=\cFin$.

\bexm
Let $\kappa$ be a cardinal number. A \emph{$\kappa$-scale in $\NN$} is a $\les$-unbounded set $\sset{a_\alpha}{\alpha<\kappa}$ in $\roth$ 
such that $a_\alpha\les a_\beta$ for all ordinal numbers $\alpha<\beta<\kappa$.
Blass~\cite{BlassMO} points out that, by a theorem of Hechler~\cite{hechler},
it is consistent that, for example,
$\fb=\aleph_1$ and there is an $\aleph_2$-scale $X=\sset{a_\alpha}{\alpha<\aleph_2}$.
Since the cardinal number $\aleph_2$ is regular,
every subset of the $\aleph_2$-scale $X$ of cardinality $\fb$ is $\les$-bounded (indeed, by some member of the set $X$). In particular, 
the set $X$ does not contain a $\fb$-unbounded set.
\eexm

\subsection{Countably Menger and countably Hurewicz spaces}

Let $\fb:=\bof(\cFin)$ and $\fd:=\bof(\roth)$. 
Applying the results of the earlier subsections to the filter $\cFin$ and the superfilter $\roth=\cFin^+$, 
we obtain the following results.

\bthm\label{thm:d-conxHur=ctblMen}\mbox{}
\be
\itm Every $\fd$-concentrated Tychonoff space satisfies
$\productive{\Hur}{\ctbl\Men}$.
\itm $\productive{\fb\con}{\fb\lind}\sub\ctbl\Men$.	
\itm $\productive{\fd\con}{\fd\lind}\sub\ctbl \Men$. In particular,
$\productive{\ctbl\Men}{\fd\lind}\sub\ctbl\Men$.
\itm If $\fb=\fd$, then
$\productive{\fd\con}{\fd\lind}\sub \ctbl\Hur$ and, in particular, 
$\productive{\ctbl\Men}{\fd\lind}\sub \ctbl\Hur$.	
\ee
\ethm
\bpf
(1) Apply Theorem~\ref{thm:d-concxH=M}.

(2) Apply Theorem~\ref{thm:ctblF2}(1).

(3) Apply Theorem~\ref{cor:ctblF+}. It is a simple, folklore observation that
every $\fd$-concentrated space is countably Menger.

(4) Apply any of the theorems~\ref{thm:ctblF}(2) or~\ref{thm:ctblS}.
\epf
	
\section{General covers}

\subsection{General semifilters}

A space is \emph{concentrated} if it is $\aleph_1$-concentrated.
Theorem~\ref{thm:ctblS} and Proposition~\ref{prp:LindctblSMen=SMen} 
imply the following result.

\bthm\label{thm:genS}
Let $S$ be a semifilter with a cofinal $S$-scale, and $\bof(S)=\aleph_1$.
\[
\productive{\Con}{\Lind}\sub S\men.\qed
\]	
\ethm

\bthm\label{thm:prodS}
Assume that $\fb=\aleph_1$.
\be
\itm Let $S$ be a meager semifilter. 
Every $S$-scale space is productively Hurewicz.
\itm Every cofinal $\cFin$-scale space is productively $S$-Menger,
for all semifilters $S$.
\ee
\ethm
\bpf 
(1) Let $X\cup \Fin$ be an $S$-scale space, and $Y$ be a Hurewicz space. By Proposition~\ref{prp:ctblProdS}(1), the product space $(X\cup\Fin)\x Y$ is countably Hurewicz. 

\blem\label{lem:ConcxH=Lind}
Every Tychonoff concentrated space satisfies $\productive{\Hur}{\Lind}$.
\elem
\bpf
This is a routine application of the projection method~\cite{AddGN}, 
see the proof of Theorem~\ref{thm:d-concxH=M}.
\epf

As the semifilter $S$ is meager, we have $\bof(S)=\fb=\aleph_1$~\cite[Corollary~2.27]{sfh}.
Being an $S$-scale, the set $X$ is $\bof(S)$-unbounded. By Lemma~\ref{lem:MichaelConc}, the space $X\cup \Fin$ is concentrated.
Therefore, the space $(X\cup\Fin)\x Y$ is also Lindel\"of, and thus Hurewicz.

(2) Let $X\cup \Fin$ be a cofinal $\cFin$-scale space, and $Y$ be an $S$-Menger space. By Proposition~\ref{prp:ctblProdS}(2), the product space $(X\cup\Fin)\x Y$ is countably $S$-Menger.

\blem\label{lem:CssxFMen=Lind}
Let $S$ be a semifilter, and $X$ be a subset of $\roth$ such that every $\le_S$-bounded subset of the set $X$ is countable. 
The space $X\cup \Fin$, with the Michael topology,
satisfies $\productive{S\men}{\Lind}$.
\elem

\bpf
Let $\le_\bbR$ be the usual order on the real line, and \[
C:=(\PN\x\{0\})\cup (\roth\x\{-1,1\})
\]
be a space with the order topology generated by the lexicographic order on $C$. 
The space $C$ is the Dedekind compactification of the space $\PN$ with the Michael topology.
Since the subspace $(X\cup\Fin)\x\{0\}$ of the space $C$ is homeomorphic to the space
$X\cup\Fin$, it suffices to prove that the space
$(X\cup\Fin)\x\{0\}$ satisfies $\productive{S\men}{\Lind}$.

Let $Y$ be an $S$-Menger space, and $\cU$ be an open cover of the product space 
\[
((X\cup\Fin)\x\{0\})\x Y,
\]
by sets open in the product space $C\x Y$.
The product space $(\Fin\x\{0\})\x Y$ is Lindel\"of. Thus, there is a countable family $\cU'\sub\cU$ 
such that the set $U:=\Un\cU'$ contains $(\Fin\x\{0\})\x Y$.  Since the set $(C\x Y)\sm U$ is closed 
in the space $C\x Y$, it is $S$-Menger. 
Thus, the projection $H$ of the set $(C\x Y)\sm U$ in $\roth\x\{-1,0,1\}$ is $S$-Menger, too.

The map $f\colon C\to\roth$ such that $f(x,y)=x$ for all elements $(x,y)\in C$, is continuous. 
Thus, the set $f[H]$ is $S$-Menger. Therefore, there is a $\le_S$-bound  $b$ for the set $f[H]$ in $\roth$.
Since the set  $\sset{x\in X}{x\le_S b}$ is countable, 
the set $X\cap f[H]$ is countable, and the set $(X\x\{0\})\cap H$ is countable, too. 
Thus, the product space $((X\x\{0\})\cap H)\x Y$ is Lindel\"of. 
Let $\cU''\sub\cU$ be a countable cover of $((X\x\{0\})\cap H)\x Y$. 
The family $\cU'\cup\cU''$ is a countable cover of the space 
$((X\cup \Fin)\x\{0\})\x Y$.
\epf

Since $\fd=\aleph_1$, every uncountable subset of the cofinal $\cFin$-scale $X$ is dominating in $\roth$. Thus, each $\leS$-bounded subset of $X$ is countable. 
Apply Lemma~\ref{lem:CssxFMen=Lind} and Proposition~\ref{prp:LindctblSMen=SMen}.
\epf

Concentration is necessary in Lemma~\ref{lem:ConcxH=Lind},
since the product of an 
uncountable space and a Hurewicz space need not be Lindel\"of. Indeed, an
uncountable discrete space is not Lindel\"of.

Theorem~\ref{thm:prodS}(2) is a generalization of an earlier result~\cite[Theorem~6.2]{SPMProd}, from hereditarily Lindel\"of spaces to general spaces.

Assuming \CH{}, every productively Lindel\"of metric space is $\sigma$-com\-pact~\cite{alster}. By Theorem~\ref{thm:prodS}(2) and Remark~\ref{rem:coarser}, we have the following corollary.

\bcor
Assume \CH{}. Let $S$ be a semifilter, and $X$ be a cofinal $\cFin$-scale in $\roth$. The real space $X\cup\Fin$ is productively $S$-Menger, but not productively Lindel\"of.\qed
\ecor

\subsection{Superfilters}
We have the following corollary of theorems~\ref{thm:ctblSup} and~\ref{thm:genS}.

\bcor\label{cor:genSup}
Let $S$ be a superfilter with $\bof(S)=\aleph_1$.
\be
\itm $\productive{\text{(cofinal) } S\scsp}{\Lind}\sub S\men$.
\itm $\productive{\Con}{\Lind}\sub S\men$.\qed
\ee
\ecor

\subsection{Filters}

By Theorem~\ref{thm:ctblF}(2) and Proposition~\ref{prp:LindctblSMen=SMen}, we have the following theorem.

\bthm\label{thm:genF}
Let $F$ be a filter such that $\bof(F^+)=\aleph_1$. 
\be
\itm  $\productive{\text{cofinal }F\scsp}{\Lind}\sub F^+\men$.
\itm  $\productive{\text{(cofinal) } F^+\scsp}{\Lind}\sub F\men$.\qed	
\ee
\ethm

\bthm\label{thm:ProdF}
Let $F$ be a filter such that $\bof(F)=\aleph_1$. 
Every cofinal $F$-scale space is productively $F$-Menger.
\ethm

\bpf
Let $X\cup\Fin$ be a cofinal $F$-scale space, and $Y$ be an $F$-Menger space. By Proposition~\ref{prp:ctblProdF}(2), the product space $(X\cup\Fin)\x Y$ is countably $F$-Menger.
Apply Lemma~\ref{lem:CssxFMen=Lind} and Proposition~\ref{prp:LindctblSMen=SMen}.
\epf

\bprp Let $U$ be an ultrafilter such that $\bof(U)=\aleph_1$.
\be
\itm $\productive{\text{(cofinal) }U\scsp}{\Lind}\sub U\men$.
\itm $\productive{\Con}{\Lind}\sub U\men$.
\itm $\productive{ U\men}{\Lind}\sub U\men$.
\itm  Every $U$-scale space is productively $U$-Menger.
\ee
\eprp
\bpf
(1) Apply Corollary~\ref{cor:genSup}(1) or Theorem~\ref{thm:genF}.

(2) Apply one of the corollaries~\ref{cor:ctblU}(1) or~\ref{cor:genSup}. 

(3) Apply (1).

(4) Apply Theorem~\ref{thm:ProdF}.
\epf

\subsection{Hurewicz, Menger, and Lindel\"of spaces}
Let $\kappa$ be a cardinal number. A real space of cardinality at least $\kappa$ is 
\emph{$\kappa$-Luzin} if the cardinalities of its intersections with 
meager sets are all smaller than~$\kappa$. A \emph{Luzin set} is an $\aleph_1$-Luzin set. 
If $\kappa\in\{\cf(\fd),\fd\}$, then for every $\kappa$-Luzin set $L$, there is a $\fd$-concentrated 
real space $Y$ such that the product space $L\x Y$ is not Menger~\cite[Corollary~2.11]{pMReal}. 
It is unknown whether, for every Luzin set $L$, there is a Menger real space $Y$ such that the product 
space $L\x Y$ is not Menger. 
We obtain a positive resolution for \emph{general} spaces.

\bthm Let $\kappa$ be an uncountable cardinal number.
\be
\itm For every $\kappa$-Luzin set $L$, there is a $\kappa$-concentrated space $Y$ such that the product space $L\x Y$ is not $\kappa$-Lindel\"of.
\itm No Luzin space is productively Menger.
\ee
\ethm
\bpf
(1) Apply Proposition~\ref{prp:not(k-con,k-lind)}.

(2) Apply (1).
\epf

We apply the results of the previous sections to Menger and Hurewicz spaces.

\bprp\label{prp:ConcxH=Men}
Every Tychonoff concentrated space satisfies $\productive{\Hur}{\Men}$.
\eprp
\bpf
Apply Theorem~\ref{thm:d-conxHur=ctblMen}(1), Lemma~\ref{lem:ConcxH=Lind}, and Proposition~\ref{prp:LindctblSMen=SMen}.
\epf

\bprp\label{prp:MH}
Assume that $\fd=\aleph_1$. 
\be
\itm $\productive{\text{(cofinal) }\roth\scsp}{\Lind}\sub \Hur$,
\itm $\productive{\Con}{\Lind}\sub \Hur$,
\itm $\productive{\text{cofinal }\cFin\scsp}{\Lind}\sub \Men$,
\itm $\productive{\Hur}{\Lind}\sub \Men$.
\ee
\eprp
\bpf
(1) Apply Theorem~\ref{thm:genF}(2).

(2) Apply~(1) or Theorem~\ref{thm:genS}.

(3) Apply Theorem~\ref{thm:genF}(1).

(4) Apply~(3) and Theorem~\ref{thm:prodS}(1).
\epf

Assume, for this paragraph, that $\fd=\aleph_1$.
Aurichi and Tall~\cite[Theorem~23]{AurTall} improved earlier results by 
proving that every productively Lindel\"of space is Hurewicz (later, 
Tall~\cite[Section~3]{TallpD} and Repov\v{s} and Zdomskyy~\cite[Theorems~1.1 and~1.2]{RZ} proved the same result using  weaker hypotheses).
It was later shown that, \emph{in the realm of hereditarily Lindel\"of spaces},
every productively Lindel\"of space is productively Hurewicz and productively Menger~\cite[Theorem~8.2]{SPMProd}. 
In our previous paper, we proved that \emph{in that realm}, 
every productively Menger space is productively Hurewicz~\cite[Theorem~4.8]{pMReal}.
We obtain an improved result in the general realm.

\bthm\label{thm:pMispH} 
Assume that $\fd=\aleph_1$. 
\[
\Lind^\x\sub\productive{\Men}{\Lind}=\Men^\x \sub \Hur^\x.
\]
\ethm
\bpf
$\productive{\Men}{\Lind}\sub\Men^\x$: 
Let $X$ be a space. Assume that there is a Menger space $M$ such that the product space $X\x M$ is not Menger. 
By Proposition~\ref{prp:MH}(3), there is a cofinal $\cFin$-scale space $Y$ such
that the product space $(X\x M)\x Y$ is not Lindel\"of. 
By Theorem~\ref{thm:prodS}, the space  $M\x Y$ is Menger. 
In summary, the product of the space 
$X$ and the Menger space $M\x Y$ is not Lindel\"of.

$\productive{\Men}{\Lind}\sub \Hur^\x$:
Let $X$ be a space. Assume that there is a Hurewicz space $H$ such that the product space $X\x H$ is not Hurewicz. 
By Proposition~\ref{prp:MH}(2), there is a concentrated space $Y$ such that the product space $(X\x H)\x Y$ is not 
Lindel\"of. By Proposition~\ref{prp:ConcxH=Men}, the space  $H\x Y$ is Menger. 
Thus, the product of the space 
$X$ and the Menger space $H\x Y$ is not Lindel\"of.
\epf

In the realm of hereditarily Lindel\"of spaces, if $\fb=\fd$, then every 
productively Menger real space is productively Hurewicz~\cite[Theorem~4.8]{pMReal}. It is unknown whether, when $\fb=\fd$,
these classes provably coincide for real spaces~\cite[Problem~6.9]{pMReal}.

\bprp\label{prp:pHisnotpM}$\ $
\be
\itm Assume that $\fb=\aleph_1<\fd$. There is a hereditarily Lindel\"of productively Hurewicz space that 
is not productively Menger.
\itm It is consistent with \CH{} that
there is a productively Hurewicz space that is not productively Menger.
\ee
\eprp
\bpf
(1) We view $\PN$ as a subset of the real line. Let $\le_\bbR$ 
be the usual order on the real line. 
For points $a,b\in\PN$, 
let $[a,b):=\sset{x\in\PN}{a\le_\bbR x <_\bbR b}$, and $(a,b]:=\sset{x\in\PN}{a<_\bbR x \le_\bbR b}$.
The \emph{Sorgenfrey topology} (\emph{Sorgenfrey* topology})~\cite{srg} on the set $\PN$ is the topology generated by the sets $[a,b)$ ($(a,b]$), for 
$a,b\in \PN$.

Let $X\cup\Fin$ be a $\cFin$-scale space such that the open neighborhoods of the points from the set $X$ are as in the space $X\cup\Fin$ with the Sorgenfrey topology.
Let $Y$ be a subset of $X$ with $\card{Y}=\aleph_1$. Equip $Y$
with the Sorgenfrey$^*$ topology.
Since the space $Y$ is Lindel\"of and $\card{Y}<\fd$, it is Menger. 
The diagonal set $\sset{(y,y)}{y\in Y}$ is a closed and discrete subset of the product space $(X\cup\Fin)\x Y$. The space $(X\cup\Fin)\x Y$ is not Lindel\"of, and thus not Menger.
By Theorem~\ref{thm:prodS}(1) and Remark~\ref{rem:coarser}, the space $X\cup\Fin$ is productively Hurewicz.

(2) We use the following lemma, pointed out to us by A. Miller.
We sketch a proof that assumes familiarity with the method of forcing.

\blem
\label{lem:arnie}
It is consistent with \CH{} that some $\cFin$-scale is Menger.
\elem
\bpf
Let $V$ be a model of \CH{}. Let $X$ be a $\cFin$-scale in $V$.
Extend $V$ generically by adding $\aleph_1$ Cohen reals.
Every Borel map $X\to\NN$ is coded in an intermediate extension,
and is thus $\lei$-bounded by the next added Cohen reals.
It follows that, in the extension, the set $X$ is a Menger space,
even with respect to countable Borel covers.
\epf

Let $X\sub\roth$ be a $\cFin$-scale with Menger's property.
Equip the space $Y:=X\cup\Fin$ with the Michael topology. 
By Theorem~\ref{thm:prodS}(1) and Remark~\ref{rem:coarser}, the space $Y$ is productively Hurewicz.
The product space $Y\x X$ contains the closed discrete subset 
$\sset{(x,x)}{x\in X}$, and is thus not Lindel\"of. In particular, 
the product is not Menger. 
\epf

\bcor\label{cor:pMisnotpH}
It is consistent with \CH{} that all inclusions in Theorem~\ref{thm:pMispH} are strict. \qed
\ecor

\section{Countably compact spaces}

A space is \emph{countably compact} if every countable open cover of this space has a finite subcover.
Let $S$ be a semifilter. Countably compact spaces are countably $S$-Menger.
Every compact space is productively $S$-Menger and productively countably $S$-Menger. 
In contrast, a result of Frolik~\cite[\S3.1.6]{Frolik}, applied to 
a discrete space of cardinality $\fd$, implies that 
there are two countably compact spaces whose product is not countably Menger.

A \emph{P-space} is a space where all $G_\delta$ sets are open.
Using results of Galvin~\cite[page~157]{GN} and Alster~\cite[Theorem~1]{alster}, 
Babinkostova, Scheepers, and Pansera~\cite[Lemma~18, Theorem~23]{WCP}
proved that every Lindel\"of P-space is productively Lindel\"of, 
productively Megner, and productively Hurewicz. 

\bprp
\label{prp:cc}
Assume that $\fd=\aleph_1$. There are a Lindel\"of P-space $X$
and a first countable, countably compact space $Y$ such that the product space
$X\x Y$ is not countably Menger.
\eprp
\bpf
Let $X:=\aleph_1+1$, where the neighborhoods of the point $\aleph_1\in X$ 
are those of the ordinal topology on the ordinal number $\aleph_1+1$,
and the remaining points are isolated.
The space $X$ is a Lindel\"of P-space.
The space $Y:=\aleph_1$, with the ordinal topology, is first countable 
and countably compact.

The diagonal set $\sset{(\alpha,\alpha)}{\alpha<\aleph_1}$ 
is a closed discrete subset of the product space $X\x Y$.
Since $\fd=\aleph_1$, the diagonal set is not countably Menger.
\epf

\bdfn
Let $X$ and $Y$ be spaces.
A set-valued map $\Phi\colon X\Rightarrow Y$ is \emph{countably compact-valued upper semicontinuous} (\emph{ccusco}) if it satisfies the definition of cusco map,
with \emph{compact} replaced by \emph{countably compact}.
\edfn

Using arguments similar to those used for cusco maps~\cite[Lemma~1]{SF1}, 
we obtain the following result.
 
\bprp\label{prp:ccusco}
Let $S$ be a semifilter.
The class of countably $S$-Menger spaces is 
preserved by ccusco maps.\qed
\eprp

Let $S$ be a semifilter. Let $\add(S\men)$ be the minimal number of $S$-Menger subspaces of $\roth$ whose union is not $S$-Menger. For a filter $F$, we have  $\add(F\men)=\bof(F)$.

\blem\label{lem:add}
Let $S$ be a semifilter. Every space that is a union of less than $\add(S\men)$ countably $S$-Menger spaces is countably $S$-Menger.\qed
\elem

A space is \emph{sequential} if every subset that is closed under limits of
sequences is closed.
In particular, every first countable space is sequential. 
By Proposition~\ref{prp:cc}, the product of countably compact spaces
with countably Menger spaces need not be countably Menger.
Item~(2) of the following results implies that 
every product of countably compact spaces
with countably Menger \emph{sequential} spaces is countably Menger.

\bprp
Let $S$ be a semifilter.
\be
\itm Every space that is a union of less than $\add(S\men)$  compact spaces is productively countably $S$-Menger.
\itm Every product of a union of less than $\add(S\men)$ countably compact spaces and a countably $S$-Menger sequential space is countably $S$-Menger.
\ee
\eprp
\bpf
(1) Apply Proposition~\ref{prp:ccusco} and Lemma~\ref{lem:add}.

(2) Let $X$ be a countably compact space, and $Y$ be a countably $S$-Menger sequential space. 
The projection $p\colon X\x Y\to Y$ is a closed map onto $Y$~\cite[Theorem~3.10.7]{eng}. 
Thus, the inverse map $p\inv\colon Y\Rightarrow X\x Y$ is ccusco.
By Proposition~\ref{prp:ccusco}, the products space $X\x Y$ is
countably $S$-Menger.
Apply Lemma~\ref{lem:add}. 
\epf

\section{Comments and open problems}

Assume \CH{}. Consider all potentially provable assertions 
\[
\productive{\text{\bf A}}{\text{\bf B}}\sub \text{\bf C},
\]
for $\text{\bf A}, \text{\bf B}, \text{\bf C}\in \{\Con, \Lind, \Men, \Hur\}$.
Some of these assertions fail for obvious reasons. 
Some others hold for trivial reasons. 
We also have the following observation.

\bprp\mbox{}
\be
\itm Assuming \CH{}, we have $\productive{\Hur}{\Men}\nsubseteq\Hur$.
\itm $\productive{\Con}{\Men}\nsubseteq\Con$.
\ee
\eprp
\bpf
(1) By Proposition~\ref{prp:ConcxH=Men}, every Luzin set satisfies $\productive{\Hur}{\Men}$. No Luzin set is Hurewicz~\cite[page~196, footnote~1]{Hure27}.

(2) Consider an uncountable non-concentrated compact space.
\epf

This shows that all provable results of the above-considered type are included
in items~(2) and~(4) of Proposition~\ref{prp:MH}.

The following problem is motivated by Theorem~\ref{thm:pMispH}.

\bprb
Assume \CH{}. 
\be
\itm Does $\productive{\Hur}{\Lind}\sub \Men^\x$?
\itm And if not, does $\productive{\Hur}{\Lind}\sub \Hur^\x$?
\ee
\eprb

A \emph{scale} is a $\les$-increasing sequence $\set{s_\alpha}{\alpha<\fd}\sub\roth$ that is $\lei$-unbounded in
$\roth$.
If $S$ is a scale then, in the realm of hereditarily Lindel\"of spaces,
the real space $S\cup\Fin$ is productively Hurewicz and productively Menger~\cite[Theorem~6.2]{SPMProd}.
A positive solution of the following problem implies the same assertion
for general spaces.

\bprb
Assume that $\fb=\fd$.
Is every real space of cardinality smaller than $\fb$
provably productively Hurewicz? Productively Menger?
\eprb

Proposition~\ref{prp:pHisnotpM} and Lemma~\ref{lem:arnie} 
motivate the following problem. A positive solution of its second item
implies a positive solution of its first item.

\bprb
\mbox{}
\be
\itm Does \CH{} imply the existence of a productively Hurewicz space that
is not productively Menger?
\itm Does \CH{} imply that some $\cFin$-scale is Menger?
\ee
\eprb

An uncountable subspace of the real line is 
\emph{Sierpi\'nski} if the cardinalities of its intersections with Lebesgue measure zero sets are all countable.
Every Sierpi\'nski set is Hurewicz~\cite[Theorem~2.10]{coc2}. Assuming \CH{}, there is a Sierpi\'nski set whose square is not Menger~\cite{coc2}. The following problem is also open in the realm of hereditarily Lindel\"of spaces~\cite[Problem~6.8]{pMReal}. 
The solution for general spaces may turn out different.

\bprb
Does \CH{} imply that no Sierpi\'nski set is productively Hurewicz?
\eprb

\subsection*{Acknowledgments}
We thank Arnold Miller for Lemma~\ref{lem:arnie}, and Franklin Tall
for pointing out some useful references.
The research of the first named author was supported by an \emph{Etiuda~2} grant,
Polish National Science Center, UMO-2014/12/T/ST1/00627.

\end{document}